\documentclass{article}
\usepackage{graphicx}
\usepackage{amsmath}
\usepackage{amssymb}
\usepackage{amscd}
\usepackage{mathrsfs}

\newtheorem{theorem}{Theorem}

\newtheorem{theoremc}{Theorem}
\newtheorem{theoremd}{Theorem}

\newtheorem{rk}[theoremc]{Remark}
\newtheorem{cor}[theoremd]{Corollary}

\newcommand\bib[1]{\bibitem[#1]{#1}}

\newcommand\abz{\hspace{13.5pt}}
\newcommand\qed{\phantom{\underline{y}}\hfill\hfill$\square$}
\newcommand{\comm}[1]{}
\newcommand\Cc{\let\mathcal\mathscr\mathcal C}

\renewcommand\a{\alpha}
\renewcommand\b{\beta}

\newcommand\C{{\mathbb C}}

\renewcommand\d{\delta}
\newcommand\D{{\mathcal D}}

\newcommand\E{{\mathcal E}}

\newcommand\g{\gamma}

\renewcommand\l{\lambda}

\newcommand\oo{\omega}
\newcommand\op[1]{\mathop{\rm #1}\nolimits}
\newcommand\ot{\otimes}
\newcommand\p{\partial}

\newcommand\R{{\mathbb R}}

\newcommand\ti{\tilde}

\newcommand\we{\wedge}

\renewcommand\v{{\rm v}}
\newcommand\z{\sigma}
\newcommand\Z{{\mathbb Z}}

\begin{document}

 \title{Differential invariants\\ of the motion group actions.}
 \author{Boris Kruglikov, Valentin Lychagin}
 \date{}
 \maketitle

 \begin{abstract}
Differential invariants of a (pseudo)group action can vary when
restricted to invariant submanifolds (differential equations). The
algebra is still governed by the Lie-Tresse theorem, but may change a lot.
We describe in details the case of the motion group $O(n)\ltimes\R^n$ acting on the
full (unconstraint) jet-space as well as on some invariant equations.%
 \end{abstract}

\section*{Introduction}

Let $G$ be a pseudogroup acting on a manifold $M$ or a bundle
$\pi:E\to M$. This action can be prolonged to the higher jet-spaces $J^k(\pi)$
(one can also start with an action in some PDE system $\E\subset J^k(\pi)$ and prolong it).

The natural projection $\pi_{k,k-1}:J^k(\pi)\to J^{k-1}(\pi)$ maps
the orbits in the former space to the orbits in the latter. If the
pseudogroup is of finite type (i.e. a Lie group), this bundle
(restricted to orbits) is occasionally a covering outside the
singularity set. Otherwise it will become a sequence of bundles for $k\gg1$.
Ranks of these bundles varies but it is occasionally given by the
Hilbert-Poincar\'e polynomial of the pseudogroup action.

The orbits can be described via differential invariants, i.e.
invariants of the action on some jet level $k$. Existence and
stability of the above mentioned Hilbert-Poincar\'e polynomial is
a consequence of the Lie-Tresse theorem, which claims that the
algebra of differential invariants is finitely generated via the
algebraic-functional operations and invariant derivations.

This theorem in the ascending degree of generality was proved in
different sources \cite{Lie$_1$,Tr,O,Ku,KL$_1$}. In particular,
the latter reference contains the
full generality statement, when the pseudogroup acts on a
system of differential equations $\E\subset J^l(\pi)$ (the standard regularity
assumption is imposed, which is an open condition in finite jets).

In the case the pseudogroup $G$ acts on the jet space, $\E$ must be
invariant and so consist of the orbits, or equivalently it has an
invariant representation $\E=\{J_1=0,\dots,J_r=0\}$, where $J_l$ are
(relative) differential invariants. Now the following dichotomy is possible.

If the orbits forming $\E$ are regular, the structure of the algebra
of differential invariants on $\E$ can be read off from that one of
the pure jet-space.

On the other hand if $\E$ consists of singular orbits\footnote{In
this case $\E$ can be defined via vanishing of an invariant tensor
${\bf J}$, with components $J_i$, though in general the latter
cannot be chosen as scalar differential invariants.} (which is often
the case when the system is overdetermined, so that differential
syzygy should be calculated, which is an invariant count of
compatibility conditions), then the structure of the algebra of
differential invariants is essentially invisible from the
corresponding algebra $\mathcal{I}$ of the pure jet-space, because
$\E$ is the singular locus for differential invariants
$I\in\mathcal{I}$ (if these exist, cf. just remarked).

In this note we demonstrate this effect on the example of motion
group $G$ acting naturally on the Euclidean space $\R^n$. The group is finite
dimensional, but even in this case the described effect is visible.
For infinite pseudogroups this follow the same route (see, for
instance, the pseudogroup of all local diffeomorphisms acting on the
bundle of Riemannian metrics in \cite{K}).

We lift the action of $G$ to the jets of functions on $\R^n$ and describe
in details the structure of algebra of scalar differential invariants in the
unconstrained ($J^\infty\R^n$) and constrained (system of PDEs) cases.
This motion group was a classical object of investigations (see e.g. the
foundational work \cite{Lie$_2$}), but we have never seen the complete
description of the differential invariants algebra.

\section{Differential invariants and Lie-Tresse theorem}

We refer to the basics on pseudogroup actions to \cite{Ku,KL$_2$},
but recall the relevant theory about differential invariants
(see also \cite{Tr,O,KJ}). Since we'll be
concerned with a Lie group in this paper, it will be denoted by one
symbol $G$ (in infinite case $G$ should be co-filtered as the equations in formal
theory).

A function $I\in C^\infty(J^\infty\pi)$ (this means that $I$ is a
function on a finite jet space $J^k\pi$ for some $k>1$) is called a differential
invariant if it is constant along the orbits of the lift of the
action of $G$ to $J^k\pi$. For connected groups $G$ we have an
equivalent formulation: The Lie derivative vanishes $L_{\hat X}(I)=0$ for
all vector fields $X$ from the lifted action of the Lie algebra.

Note that often functions $I$ are defined only locally near families
of orbits. Alternatively we should allow $I$ to have meromorphic
behavior over smooth functions (but we'll be writing though about
local functions in what follows, which is a kind of micro-locality,
i.e. locality in finite jet-spaces).

The space $\mathcal{I}=\{I\}$ forms an algebra with respect to usual
algebraic operations of linear combinations over $\R$ and
multiplication and also the composition $I_1,\dots,I_s\mapsto
I=F(I_1,\dots,I_s)$ for any $F\in C^\infty_\text{loc}(\R^s,\R)$,
$s=1,2,\dots$ any finite number. However even with these operations
the algebra $\mathcal{I}$ is usually not locally finitely generated. Indeed,
the subalgebras $\mathcal{I}_k\subset\mathcal{I}$ of order $k$
differential invariants are finitely generated on non-singular
strata with respect to the above operations, but their injective
limit $\mathcal{I}$ is not.

To cure this difficulty S.Lie and later his French student A.Tresse
introduced invariant derivatives, i.e. such differentiations
$\vartheta$ that belong to the centralizer of the Lie algebra
$\mathfrak{g}=\op{Lie}(G)$ lifted as the space of vector fields on
$J^\infty(\pi)$. To be more precise we consider the derivations
$\vartheta\in C^\infty(J^\infty\pi)\ot_{C^\infty(M)}\mathcal{D}(M)$
($\Cc$-vector fields on $\pi$), which commute with the $G$-action.
These operators map differential invariants to differential invariants
$\vartheta:\mathcal{I}_k\to\mathcal{I}_{k+1}$.

We can associate invariant differentiations to a collection of
differential invariants $I_1,\dots,I_n$ ($n=\dim M$) in general
position, meaning $\hat d I_1\we\dots \we\hat d I_n\ne0$. Moreover
the whole theory discussed above transforms to the action on
equations\footnote{At this point we do not need to require even formal
integrability of the system $\E$ \cite{KL$_1$}, but this as well as regularity
issues will not be discussed here.} $\E\subset J^\infty(\pi)$.

Namely, given $n$ functionally independent invariants $I^1,\dots,I^n$ we
assume their restrictions $I^1_\E,\dots,I^n_\E$ are functionally
independent\footnote{Here and in what follows one can assume (higher micro-)local
treatment.} (in fact we can have the latter invariants only without the former),
so that they can be considered as local coordinates.

Then one can introduce the horizontal basic forms (coframe)
$\oo^i=\hat d I^i_\E$. Its dual frame consists of invariant
differentiations $\hat\p/\hat\p
I^i_\E=\sum_j[\D_a(I^b_\E)]^{-1}_{ij}\D_j$. The invariant derivative
of a differential invariant $I$ are just the coefficients of the
decomposition of the horizontal differential by the coframe:
 $$
\hat d I=\sum_{i=1}^n\frac{\hat\p I}{\hat\p I^i_\E}\,\oo^i
 $$
and they are called Tresse derivatives.

All invariant tensors and operators can be expressed through the given frame and
coframe and this is the base for the solution of the equivalence
problem.

Lie-Tresse theorem claims that the algebra of differential
invariants $\mathcal{I}$ is finitely generated with respect to
algebraic-functional operations and invariant derivatives.

\section{Motion group action}

Consider the motion group $\op{O}(n)\ltimes\R^n$. It is disconnected
and for the purposes of further study of differential invariants we
restrict to the component of unity $G=\op{SO}(n)\ltimes\R^n$. The
two Lie groups have the same Lie algebra
$\mathfrak{g}=o(n)\ltimes\R^n$ and the differential invariants of
the latter become the differential invariants of the second via
squaring.

Since the latter is inevitable even for the group $G$, the difference
between two algebras of invariants is by an extension via finite group and
will be ignored.

Below we will make use of the action of $G$ on the space
of codimension $m$ affine subspaces of $\R^n$:
 $$
\op{AGr}(m,n)\equiv\{\Pi+c\}\simeq\{(\Pi,c):\Pi\in\op{Gr}(n-m,n),c\in\Pi^\perp\}.
 $$
The action of $G$ is $x\mapsto Ax+b$, $x\in\R^n$, it is transitive
on $\op{AGr}(m,n)$ and the stabilizer equals
 $$
\op{St}(\Pi+c)=\{(A,b)\in G:A\Pi=\Pi,b\in(1-A)c+\Pi\}\simeq
\op{SO}(\Pi)\times\op{SO}(\Pi^\perp)\ltimes\Pi.
 $$
We have $\dim G=\dfrac{n(n+1)}2$, $\dim\op{AGr}(m,n)=m(n-m+1)$ and
 $$
\op{AGr}(k,n)\simeq G/(\op{SO}(m)\times\op{SO}(n-m)\ltimes\R^{n-m})
 $$
(note that this implies $\op{AGr}(m,n)\ne\op{AGr}(n-m,n)$ except for
$n=2m$ contrary to the space $\op{Gr}(m,n)$).

We can extend the action of $G$ on $\R^n$ to the space
$\R^n\times\R^m$ by letting $g\in G$ act
 $$
g\cdot(x,u)=(g\cdot x,u).
 $$
We can prolong the action to the space $J^k(n,m)$.

For $k=1$ the action commutes with the natural $\op{Gl}(m)$-action in
fibers of the bundle $\pi_{10}:J^1(n,m)\to J^0(n,m)$ and the action
descends on the projectivization, which can be identified with the
open subset in $\R^n\times\op{AGr}(k,n)$ by associating the space
$\op{Ker}(d_xf)$ to a (surjective at $x$ if we assume $n>m$)
function $f:\R^n\to\R^m$.

Thus $u$ is indeed an invariant of the $G$-action (scalar invariants
are its components $u^i$, so that we can assume the fiber $\R^m$
being equipped with coordinates), and the scalar differential
invariants of order 1 are\footnote{Recall that the base space $\R^n$
is equipped with the Euclidean metric preserved by $G$.}
$\langle\nabla u^i,\nabla u^j\rangle=\sum u^i_{x^s}u^j_{x^s}$.

These form the generators of scalar differential invariants of
order\footnote{This claim holds at an open
dense subset of $J^1(n,m)$. However if we restrict to the set of
singular orbits with $\op{rank}(d_xu)=r<m$, the basic set of invariants
will be quite different.} $\le1$.

 \begin{rk}
Sophus Lie investigated the vertical actions of $G$ in $J^0(m,n)=\R^m\times\R^n$
and the invariants of its lift to $J^\infty(m,n)$ \cite{Lie$_2$} (actually in this paper
for $m=1,n=3$). This case is easier
since the total derivatives $\D_1,\dots,\D_m$ are obvious invariant derivations.
 \end{rk}

In what follows we restrict to the case $m=1$ and investigate
invariants of the $G$-action in $J^\infty(n,1)=J^\infty(\R^n)$.
Partially the results
extend to the case of general $m$, though the theory of
vector-valued symmetric forms $S^k(\R^n)^*\ot\R^m$ is more
complicated.

\section{Differential invariants: Space $J^\infty(\R^n)$}

Denote $V=T_0\R^n$. Our affine space $\R^n$ (as well as the vector space $V$)
is equipped with the Euclidean scalar product $\langle,\rangle$
and $G$ is the symmetry group of it. In
what follows we will identify the tangent space $T_x\R^n$ with $V$
via translations (using the affine structure on $\R^n$).

The space $J^\infty(\R^n)$, which is the projective limit of the
finite-dimensional manifolds $J^k(\R^n)$, has coordinates
$(x^i,u,p_\z)$, where $\z=(i_1,\dots,i_n)\in\Z_{\ge0}^n$ is a
multiindex with length $|\z|=i_1+\dots+i_n$.

The only scalar differential invariants\footnote{From now
on by this we mean the minimal set of generators.} of order $\le1$
are
 $$
I_0=u \text{ and } I_1=|\nabla u|^2.
 $$

For each $x_1\in J^1(\R^n)$ the group $G$ has a large stabilizer.
Provided $x_1$ is non-singular the dimension of the stabilizer $\op{St}_1$ is
$\dim G-2n+1=\frac12(n-1)(n-2)$.

However the stabilizer completely evolves upon the next
prolongation: the action of $G$ on an open dense subset of $J^k(\R^n)$
for any $k\ge2$ is free. Note that due to the trivial connection in
$J^0(\R^n)=\R^n\times\R$ we can decompose
 \begin{equation}\label{split}
J^k(\R^n)=\R^n\times\R\times V^*\times S^2V^*\times\dots\times S^kV^*.
 \end{equation}
Thus we can represent a point $x_k\in J^k(\R^n)$ as the base projection
$x\in\R^n$ and a sequence of "pure jets" $Q_t=d^tu\in S^tV^*$,
$t=0,\dots,k$.

Covector $Q_1$ can be identified with the vector ${\rm v}=\nabla u$.

Consider the quadric $Q_2\in S^2V^*$. Due to the metric we can
identify it with a linear operator $A\in V^*\ot V$, which has
spectrum
 $$\op{Sp}(A)=\{\l_1,\dots,\l_n\}$$
and the normalized eigenbasis $e_1,\dots,e_n$ (each element defined up to a sign!),
provided $Q_2$ is semi-simple. Since $Q_2$ is symmetric, the basis is orthonormal.

In what follows we assume to work over the open dense subset $U\subset
J^2(\R^n)$, where $A$ is simple, so that the basis is defined (almost) uniquely
(this can be relaxed to semi-simplicity, but then the stabilizer is
non-trivial and the number of scalar invariants drops a bit).

There are precisely $(2n-1)=\dim J^2(\R^n)-\dim\op{St}_1$ differential
invariants of order 2. One choice is to take ${\bar I}_{2,i}=\l_i$
and ${\bar I}_{2,(i)}=\langle e_i,{\rm v}\rangle$, $i=1,\dots,n$. There is an
obvious relation $\sum_{i=1}^n{\bar I}_{2,(i)}^{\,2}=1$, so that we can
restrict to the first $(n-1)$ invariants in this group, but beside
this the invariants are functionally independent.

Another choice of invariants is provided by the restriction $Q_\Pi$
of $Q_2$ to $\Pi={\rm v}^\perp$, which has spectrum (again by
converting quadric to an operator)
$\op{Sp}(Q_\Pi)=\{\ti\l_1,\dots,\ti\l_{n-1}\}$ and normalized
eigenvectors $\tilde e_i$. So the following invariants can be
chosen: ${\ti I}_{2,i}=\ti\l_i$, ${\ti I}_{2,n}=Q_2({\rm v},{\rm
v})$ and ${\ti I}_{2,(i)}=Q_2({\rm v},\ti e_i)$.

Both choices have disadvantages of using transcendental functions
(solutions to algebraic equations), but we can overcome this with
the following choice:
 $$
I_{2,i}=\op{Tr}(A^i),\ I_{2,(i)}=\langle A^i{\rm v},{\rm
v}\rangle,\quad i=1,\dots,n.
 $$
Here the number of invariants is $2n$, but they are
dependent\footnote{The first $(2n-1)$ invariants are however
independent and algebraic in the jets.} due to Newton-Girard
formulas, which relate the elementary symmetric polynomials
$E_k(A)=\sum_{\i_1<\dots<\i_k}\l_{i_1}\cdots\l_{i_k}$ and
power sums $S_k(A)=\op{Tr}(A^k)=\sum\l_i^k$ (these are $I_{2,k}$):
 $$
kE_k(A)=\sum_{i=1}^k(-1)^{i-1}S_i(A)E_{k-i}(A),
 $$
which together with $E_0=1$ gives an infinite chain of formulas
 $$
E_1=S_1,\ \ 2E_2=S_1^2-S_2,\ \ 6E_3=S_1^3-3S_1S_2+2S_3,\ \dots
 $$

Now with the help of Cayley-Hamilton formula
 $$
A^n=E_1(A)A^{n-1}-E_2(A)A^{n-2}+\dots+(-1)^nE_{n-1}(A)A-(-1)^nE_n(A)
 $$
we can express
 $$
I_{2,(n)}=E_1(A)I_{2,(n-1)}-E_2(A)I_{2,(n-2)}+\dots-(-1)^n\det A
 $$
through our invariants since $E_i(A)$ are functions of $I_{2,i}$.

 \begin{rk}\label{rk1}
We could restrict only to invariants $I_{2,(i)}$, $i=1,\dots,2n-1$.
This is helpful as we shall see. But when we restrict to singular
(from the orbits point of view) PDEs these differential invariants
may turn to be non-optimal, and this will be precisely the case in
the example we investigate.
 \end{rk}

Now there are precisely  $\binom{n+2}3=\dim S^3V^*$ differential
invariants of order 3, $\binom{n+3}4=\dim S^4V^*$ differential
invariants of order 4, \dots, $\binom{n+k-1}k=\dim S^kV^*$
differential invariants of order $k$.

The third order invariants are the following:
 $$
{\bar I}_{3,\z}=Q_3(e_i,e_j,e_l),\ \text{ where } \z=(ijl)\in
S^3\{1,\dots,n\}.
 $$
Generating invariants of orders 4 and higher are obtained from the
similar formulae, namely as the coefficients $q_\z$ of the
decomposition
 $$
Q_k=\sum_{\z=(i_1,\dots,i_k)}q_\z\oo^\z,\quad\text{ where
}\oo^\z=\oo^{i_1}\cdots\,\oo^{i_k},\ 1\le i_1\le\dots\le i_k\le n.
 $$

They are again transcendental functions. To get algebraic
expressions one can use the third order functions
 $$
I_{3,\z}=Q_3(A^i{\rm v},A^j{\rm v},A^l{\rm v}),\quad \z=(ijk) \text{
with } 1\le i\le j\le l\le n
 $$
and similar expressions for the higher order.

 \begin{theorem}
The invariants $I_{i,\z}$ with $i\le 3$ is the base of differential
invariants for the Lie group $G$ action in $J^\infty(\R^n)$ via
algebraic-functional operations and Tresse derivatives.
 \end{theorem}

This statement is an easy dimensional count\footnote{In fact for $n\le 4$
the same arguments imply that the base can formed only by the
invariants $I_{i,\z}$ with $i\le 2$.} together with examination of
independency condition. To get Tresse
derivatives $n$ invariants (for instance of order $\le2$) should be
chosen.

However this is not necessary, if one does not care about
transcendental functions. Indeed, the vector fields $e_1,\dots,e_n$
are invariant differentiations (they can be expressed through the
total derivatives $\D_1,\dots,\D_n$ with coefficients of the second
order).

 \begin{rk}
Notice that the moving frame
 $$e_1,\dots,e_n\in C^\infty(U,\pi_2^*T\R^n)$$
uniquely fixes an element $g\in G$, which transforms it to the
standard orthonormal frame at $0\in\R^n$. This leads to the
equivariant map defined on the open dense set
$\pi_{\infty,2}^{-1}(U)$:
 $$
J^\infty(\R^n)\to J^2(\R^n)\supset U\to G.
 $$
Such map is called the moving frame in the approach of Fells and
Olver \cite{FO}.
 \end{rk}

\section{Relations in the algebra $\mathcal{I}$}

Since the commutator of invariant differentiations is an invariant
differentiation, decomposition $[e_i,e_j]=\sum c_{ij}^ke_k$ yields
$\le \frac12n^2(n-1)$ (in general precisely this number) 3rd order
differential invariants $c_{ij}^k$.  The number of pure 3rd order
invariants obtained via invariant differentiations of the 2nd order
invariants is $n(2n-1)$. So since
 $$
\frac{n(n+1)(n+2)}6-n(2n-1)-\frac{n^2(n-1)}2=\frac{n(n+4)(1-n)}3\le0
 $$
we can conclude that differential invariants $I_{i,\z}$ with $i\le2$
and invariant differentiations $\{e_i\}_{i=1}^n$ generate the whole
algebra $\mathcal{I}$ on an open set $\hat U\subset J^\infty(\R^n)$.

Thus we are lead to the question on relations in this algebra. They
can be all deduced from the expressions for pure jets of $u$
 $$
Q_3=\hat\nabla Q_2,\ Q_4=\hat\nabla Q_3\ \text{ etc }
 $$
using the structural equations. Here
 $$
\hat\nabla:C^\infty(\pi_i^*S^iV^*)\to C^\infty(\pi_{i+1}^*S^{i+1}V^*)
 $$
is the symmetric covariant derivative
induced by the flat connection $\nabla$ in the trivial bundle
$J^0(\R^n)=\R^n\times\R$, $V=T\R^n$ (the map is
the composition of the horizontal differential $\hat d$ and
symmetrization).

However for the sake of algebraic formulations we change invariant
differentiations $e_i$ to the following ones:
 \begin{alignat*}{5}
v_1=&\hat{\rm v}={\rm v}\cdot\D_x=\sum u_i\D_i\\
v_2=&\hat{A{\rm v}}=A{\rm v}\cdot\D_x=\sum u_iu_{ij}\D_j\\
v_3=&\widehat{A^2{\rm v}}=A^2{\rm v}\cdot\D_x=\sum u_iu_{ij}u_{jk}\D_k\\
&\dots\qquad \dots\qquad \dots \\
v_n=&\widehat{A^{n-1}{\rm v}}=A^{n-1}{\rm v}\cdot\D_x=\sum
u_{i_1}u_{i_1i_2}\dots u_{i_{n-1}i_n}\D_{i_n}.
 \end{alignat*}

Now we are going to change the basis of differential invariants in
$\mathcal{I}_k$ to describe the relations in the simplest way.

Namely for the basis of invariants of order 2 we can take
$I_{2,(ij)}=Q_2(A^i\v,A^j\v)$, $0\le i\le j<n$. However since
$Q_2(v,w)=\langle Av,w\rangle$ and $A$ is self-adjoint we get
 $$
I_{2,(ij)}=\langle A^{i+1}\v,A^j\v\rangle=\langle
A^{i+j+1}\v,\v\rangle=I_{2,(i+j+1)},
 $$
so that the new invariants are precisely the old ones $I_{2,(i)}$,
just with the larger index range $i=1,\dots,2n-1$ (we can allow
arbitrary index $i$, but the corresponding
invariants are expressed via these ones, see Remark \ref{rk1} and before).

Basic higher order invariants are introduced in the same fashion:
 $$
I_{s,(i_1\dots i_s)}=Q_s(A^{i_1}\v,\dots,A^{i_s}\v),\quad 0\le
i_1\le\dots\le i_s<n.
 $$

Suppose now that our set of generic (regular) points $U\subset
J^2(\R^n)$ is given by not only the constraint that $\op{Sp}(A)$ is simple,
but also the claim that the $n\times n$ matrix $\|\g_{ij}\|_{0\le
i,j<n}$ with entries
$\g_{ij}=\langle A^i\v,A^j\v\rangle=I_{2,(i+j)}$ is non-degenerate.
Let
 $$
[\g^{ij}]=\begin{pmatrix} 1 & I_{2,(1)} & \cdots & I_{2,(n-1)} \\
I_{2,(1)} & I_{2,(2)} & \cdots & I_{2,(n)} \\
\vdots & \vdots & \ddots & \vdots \\
I_{2,(n-1)} & I_{2,(n)} & \cdots & I_{2,(2n-2)}
\end{pmatrix}^{-1}
 $$
be the inverse matrix. Note that all its entries are invariants. Now
 \begin{multline*}
(A^{i_0}\v\cdot\D_x)\, Q_s(A^{i_1}\v,\dots,A^{i_s}\v)
=Q_{s+1}(A^{i_0}\v,A^{i_1}\v,\dots,A^{i_s}\v)\\
+\sum_{j=1}^s Q_s(A^{i_1}\v,\dots,A^{i_j-1}\v,\theta_{i_0i_j},A^{i_j+1}\v,\dots,A^{i_s}\v),
 \end{multline*}
where $\theta_{i_0i_j}=\nabla_{A^{i_0}\hat\v}(A^{i_j}\v)$ is the vector which,
due to metric duality, is dual to the covector $\sum_{\a+\b=i_j-1}Q_3(A^{i_0}\v,A^\a\v,A^\b\cdot)$.
Thus we obtain

 \begin{theorem}
The algebra $\mathcal{I}$ is generated by the invariants $I_{s,\z}$
and invariants derivatives $v_1,\dots,v_n$, which are related by the
formulae ($s\ge2$):
 $$
\!\!
v_{i_0}\cdot I_{s,(i_1\dots i_s)}=I_{s+1,(i_0i_1\dots
i_s)}+\sum_{j=1}^s\sum_{a,b=0}^{n-1}\sum_{\a+\b=i_j-1}
\!\!
I_{s,(i_1\dots i_{j-1},a,i_{j+1}\dots i_s)}\g^{ab}
I_{3,(i_0,\a,b+\b)}.
 $$
In this case we can choose $I_{s,\z}$, $s\le 3$ and $v_i$ as the
generators.
 \end{theorem}
This representation for $\mathcal{I}$ via generators and relations
is not minimal, as clear from the first part of the section. However
the relations are algebraic, explicit and quite simple.

To explain how to achieve minimality let us again change the set of
generators (basic differential invariants). For the second order we
return to $I_{2,i},I_{2,(i)}$, $1\le i\le n$. For the third order we
add the invariants
 $$
I_{3,[ij]l}=\op{Tr}(Q_3(A^i\cdot,A^j\cdot,A^l\v)).
 $$
They can indeed be expressed algebraically through the invariants
$I_{3,(ijk)}$ together with the lower order invariants.

For higher order we have more possibilities of inventing new
invariants (which can be described via graphs of the type $(k,1)$-tree),
but they are again algebraically dependent with already known
differential invariants.

The relations are as follows ($0\le k< n$ and we show only top of
the list):
 \begin{gather*}
v_1\cdot I_0=I_1,\ v_2\cdot I_0=I_{2,(1)},\ \dots,\ v_n\cdot
I_0=I_{2,(n-1)},\\
v_1\cdot I_1=2I_{2,(1)},\ v_2\cdot I_1=2I_{2,(2)},\ \dots,\ v_n\cdot
I_1=2I_{2,(n)},\\
v_{k+1}\cdot I_{2,l}=\!\!\sum_{\a+\b=l-1}\!\!I_{3,[\a\b]k},\
v_{k+1}\cdot I_{2,(l)}=\!\!\sum_{\a+\b=l-1}\!\!I_{3,(\a\b k)}
 +2I_{2,(k+l+1)}\ \text{ etc.}
 \end{gather*}

Elaborate work with these shows that all the invariants can be
obtained from $I_0$ and structural constants ${\bar c}_{ij}^k$ of
the frame $[v_i,v_j]=\sum{\bar c}_{ij}^kv_k$.

 \begin{cor}
By shrinking $\hat U\subset J^\infty(\R^n)$ further (but leaving it
open dense) we can arrange that the algebra $\mathcal{I}$ of
differential invariants is generated only by $I_0$ and the
derivations $v_1,\dots,v_n$.
 \end{cor}

\section{Algebra of differential invariants: Equation $\E$}

Consider the PDE $\E=\{\|\nabla u\|=1\}$. By the standard arguments
it determines a cofiltered manifold in $J^\infty(\R^n)$ and we identify
$\E$ with it, so that it consists of the sequence of prolongations
$\E_k\subset J^k(\R^n)$ and projections
${\bar\pi}_{k,k-1}:\E_k\to\E_{k-1}$.

Since the prolongation of the defining equation for $\E$ to the
second jets is $Q_2(\v,\cdot)=0$ or $\v\in\op{Ker}(A)$ we conclude
that most of the invariants, introduced on the previously defined
subset $\hat U$, vanish: the equation is singular. Indeed,
$0\in\op{Sp}(A)$, so that $\det A=0$, the matrix $[\g_{ij}]$ is not
invertible etc.

In particular, $I_{2,(i)}=0$, $I_{s,(i_1\dots i_s)}=0$ if at least
one $i_t\ne0$, $v_2=\dots=v_n=0$. Thus the algebra $\mathcal{I}$
description from the previous section does not induce any
description of the algebra $\mathcal{I}_\E$ of differential
invariants of the group $G$ action on $\E$: the notion of regularity
and basic invariants are changed completely!

Again the group acts freely on the second jets. So there is 1
invariant of order 0
 $$
I_0=u,
 $$
no invariants of order 1 and $(n-1)$ invariants of order 2:
 $$
I_{2,1},\dots,I_{2,n-1}\qquad\text{ or equivalently }\quad
E_1(A),\dots,E_{n-1}(A).
 $$
The number of invariants of pure order $k>2$ coincides with the
ranks of the projections:
 $$\dim{\bar\pi}_{k,k-1}^{-1}(*)=\binom{n+k-2}k.$$

The principal axes of $Q_2$ (or normalized eigenbasis of $A$) are
now $e_1=\v,e_2,\dots,e_n$. These are still the invariant
derivations and the invariants of order $k>2$ are the
coefficients\footnote{Note that these invariants are defined up to $\pm$
and so should be squared to become genuine invariants; alternatively
certain products/ratios of them define absolute invariants.}
of the decomposition by basis in $S^k\op{Ann}(\v)\subset S^kV^*$:
 $$
Q_k|_\E=\sum_{\z=(i_1\dots i_k): i_t>1}q_\z\,\oo^\z,\quad
q_\z=Q_k(v_{i_1},\dots,v_{i_k}).
 $$

 \begin{theorem}
The invariants $I_0$, $I_{2,i}$ and $I_{3,\z}$ $\bigl(1\le i<n,\
\z=(i_1,i_2,i_3),\ i_t\ne1\bigr)$ form a base of differential
invariants of the algebra $\mathcal{I}_\E$ via algebraic-functional
operations and Tresse derivatives.
 \end{theorem}

Algebra of differential invariants can again be represented in a
simpler form via differential invariants and invariant derivatives.
If we choose $e_i$ for the latter the relations can be read off from
the algebra $\mathcal{I}$, though this again involves transcendental
functions.

Denote the Christoffel symbols of $\hat\nabla$ in the basis $e_\a$ by
$\Gamma_{ij}^k$ (these are differential invariants of order 3):
 $$
\hat\nabla_{e_i}e_j=\sum\Gamma_{ij}^ke_k\quad \Longleftrightarrow\quad
\hat\nabla_{e_i}\oo^j=-\sum\Gamma_{ik}^j\oo^k.
 $$
Notice that since the connection is torsionless, $T_\nabla=0$, these
invariants determine the structure functions
$c_{ij}^k=\Gamma_{ij}^k-\Gamma_{ji}^k$.

Let us now substitute the formulas (eigenvalues $\l_i$ can be
expressed through the invariants $I_{2,i}$, however in a
transcendental way; $\l_1=0$ corresponds to $e_1$)
 $$
Q_2=\sum_{1<i\le n}\l_i(\oo^i)^2,\quad Q_3=\sum_{1<i\le j\le k\le
n}q_{ijk}\oo^i\oo^j\oo^k
 $$
into the identity $\hat\nabla Q_2=Q_3$:
 \begin{multline*}
\hat\nabla\sum\l_i(\oo^i)^2=\sum(\hat\nabla\l_i)(\oo^i)^2+2\sum\l_i\oo^i\cdot\hat\nabla\oo^i\\
=\sum\p_{e_k}(\l_i)\oo^i\oo^i\oo^k-2\sum\l_i\Gamma_{jk}^i\oo^i\oo^j\oo^k.
 \end{multline*}
We get for $1<i\le j\le k\le n$:
 $$
q_{ijk}=\bigl(\p_{e_k}(\l_i)\d_{ij}+\p_{e_i}(\l_k)\d_{jk}-\p_{e_k}(\l_i)\d_{ik}\bigr)-2\sum_{\tau\in
S_3}\l_{\tau(i)}\Gamma_{\tau(j)\tau(k)}^{\tau(i)}
 $$
Since in addition, in general position the invariants $\l_i$ can be
expressed through the invariants $e_i\cdot I_0$ ($1<i\le
n$)\footnote{We have $e_1\cdot I_0=1$ on $\E$.}, then by adding
decomposition of the covariant derivatives by the frame into the set
of operations, we obtain the following

 \begin{cor}
By shrinking $\hat U\subset\E$ further (but leaving it open dense)
we can arrange that the algebra $\mathcal{I}_\E$ of differential
invariants is generated only by $I_0$ and the derivations
$e_1,\dots,e_n$.
 \end{cor}

\section{Algebra of differential invariants: Equation
$\ti\E$}\label{S5}

Completely new picture for the algebra of differential
invariants emerges, when we add one more invariant PDE: the system
becomes overdetermined and compatibility conditions (or differential
syzygies) come into the play.

We will study the following system\footnote{This interesting system
was communicated to the first author by Elizabeth Mansfield.}, which
comes from application to relativity \cite{C} (when Laplacian
$\Delta$ is changed to Dalambertian $\square$):
 $$\{\|u\|=1,\Delta u=f(u)\}\subset\E.$$
This equation is a non-empty submanifold in $J^2(\R^n)$, but when we
carry the prolongation-projection scheme, it becomes much smaller.

It turns out that for most functions $f(u)$ the resulting
submanifold $\ti\E$ is just empty. We are going to decompose it into
the strata
 $$
\ti\E=\Sigma_1(\ti\E)\cup\dots\cup\Sigma_n(\ti\E),
 $$
where $\Sigma_i(\ti\E)=\{x\in\ti\E:\ \#[\op{Sp}(A_{\ti\E})]=i\}$ for
the operator $A_{\ti\E}$ corresponding to the 2-jet $Q_2|_{\ti\E}$.

It is possible to show that the spectrum of $A$ on $\ti\E$ depends
on $u$ (and some constants) only. This was done in \cite{FZY} via
the Cayley-Hamilton theorem, though they used the Dalambertian
instead of the Laplace operator. In the next section we prove it for
the Laplace operator via a different approach.

More detailed investigation leads to the following claim:

 \begin{quote}
{\bf Conjecture:} The strata $\Sigma_n(\ti\E),\dots,\Sigma_3(\ti\E)$
are empty, while $\Sigma_2(\ti\E),\Sigma_1(\ti\E)$ are not and they
are finite-dimensional manifolds.
 \end{quote}

Let us indicate the idea of the proof for the stratum
$\Sigma_n(\ti\E)$ because on other strata the eigenbasis $e_i$ is
not defined (but the arguments can be modified). It turns out that
the compatibility is related to dramatic collapse of the algebra
$\mathcal{I}_{\ti\E}$ of differential invariants.

Indeed, as follows from the discussion above and the next section,
there is only one invariant $u$ of order $\le2$ for the $G$-action on $\ti\E$.
Since the coefficients of the invariant derivations have the second order,
we obtain the following statement:

 \begin{theorem}
All differential invariants of the Lie group $G$-action on the PDE
system $\ti\E$ can be obtained from the function $I_0=u$ and
invariant derivations.
 \end{theorem}

Now relations in the algebra $\mathcal{I}_{\ti\E}$ are differential
syzygies for $\ti\E$ and they boil down to a system of ODEs on
$f(u)$, which completely determines it.

The details of this program will be however realized elsewhere.

\section{Geometry of the system}\label{S5.5}

 \abz
In this section we justify the claim from \S\ref{S5} and prove that
the spectrum of the operator $A=A_{\ti\E}$, obtained from the pure
2-jet $Q_2|_{\ti\E}$ via the metric, depends on $u$ only. To do this
we reformulate the problem with nonlinear differential equations in
the geometric language from contact geometry \cite{Ly}.

The first equation $\E$ we represent as a level surface
$H=\frac12(1-\sum_{i=1}^np_i^2)=0$ in the jet-space $J^1(\R^n)$. The
second equation from $\ti\E$ can be represented as Monge-Ampere
type via $n$-form
 $$
\Omega_1=\sum_{i=1}^n dx_1\we\dots\we dx_{i-1}\we dp_i\we\dots\we
dx_n-f(u)dx_1\we\dots\we dx_n.
 $$
Namely a solution to the system is a Lagrangian submanifold
$L^n\subset\{H=0\}$ such that $\Omega_1|_{L^n}=0$. Representing
$L^n=\op{graph}[j^1(u)]$ we obtain the standard description.

The contact Hamiltonian vector field $X_H$ preserves the contact
structure and being restricted to the surface $H=0$ it coincides
with the field of Cauchy characteristic $Y_H=X_H|_{H=0}=\sum
p_i\D_{x_i}=\sum p_i\p_{x^i}+\p_u$.

Since Cauchy characteristics are always tangent to any solution, the
forms $\Omega_{1+i}=(L_{X_H})^i\Omega_1$ also vanish on any solution
of the system $\E$. We simplify them modulo the form $\Omega_1$ and
get:
 \begin{align*}
&\quad\Omega_2=L_{X_H}\Omega_1+f(u)\Omega_1\hfill\quad\hphantom{aaaaaaaaaaaaaaaaaaaaaaaaaaa}\\
&=2\sum dx_1\we\dots dp_i\we dx_{i+1}\dots dp_j\we dx_{j+1}\dots\we dx_n-(f'+f^2)\,dx_1\we\dots dx_n,\\
&\quad\Omega_3=L_{X_H}\Omega_2+(f'(u)+f^2(u))\Omega_1\\
&\qquad=3!\sum  dx_1\we\dots dp_i\we dx_{i+1}\dots dp_j\we dx_{j+1}\dots dp_k\we dx_{k+1}\dots\we dx_n\\
&\qquad\qquad-(D+f)^2(f)\,dx_1\we\dots\we dx_n,\\
&\qquad\dots\quad\dots\quad\dots\quad\dots\\
&\quad\Omega_n=n!\,dp_1\we\dots dp_n-(D+f)^{n-1}(f)\,dx_1\we\dots\we dx_n,\\
&\quad\Omega_{n+1}=-(D+f)^n(f)\,dx_1\we\dots\we dx_n,
 \end{align*}
where $D$ is the operator of differentiation by $u$ and $f$ is the
operator of multiplication by $f(u)$. Thus a necessary condition for
solvability is the following non-linear ODE:
 \begin{equation}\label{ODE}
(D+f)^{n+1}(1)=0.
 \end{equation}
This equation can be solved via conjugation $D+f=e^{-g}De^g$ with
$g(u)=\int f(u)\,du$ \cite{Ko}, which reduces the ODE to the form
$D^{n+1}e^g=0$, so that $g=\op{Log} P_n(u)$, where $P_n(u)$ is a
polynomial of degree $n$, whence\footnote{Here we can assume we are
working over $\C$, though this turns out to be inessential.}
 \begin{equation}\label{fu=}
f(u)=\sum_{i=1}^n\frac1{u-\a_i},\qquad \a_i=\op{const}.
 \end{equation}

However there are more compatibility conditions, which produce
further constraints on numbers $\a_i$. The above relations
$\Omega_i=0$ can be used to find $\op{Sp}(A)$. Namely let us rewrite
them as follows:
 \begin{multline*}
E_1(A)=\sum\l_i=f,\quad
E_2(A)=\sum_{i<j}\l_i\l_j=\tfrac1{2}(D+f)^2(1),\\
E_3(A)=\sum_{i<j<k}\l_i\l_j\l_k=\tfrac1{3!}(D+f)^3(1),\ \dots,\
E_n(A)=\l_1\cdots\l_n=\tfrac1{n!}(D+1)^n(1).
 \end{multline*}
These, due to Newton-Girard formulas, imply the equivalent
identities:
 \begin{multline*}
\quad I_{2,1}=\sum\l_i=f(u),\quad I_{2,2}=\sum\l_i^2=-f'(u),\\
I_{2,3}=\sum\l_i^3=\tfrac12f''(u),\quad
I_{2,4}=\sum\l_i^4=-\tfrac1{3!}f'''(u),\quad\dots
 \end{multline*}
In particular we get $\l_i=(u-\a_i)^{-1}$ and so
 $$
A\sim\op{Diag}\left(\frac1{u-\a_1},\dots,\frac1{u-\a_n}\right).
 $$
The fact that $\op{det}(A)=0$ on $\ti\E$ implies that $\a_n=\infty$
and using symmetry $\p_u$ (shift along $u$) we can arrange $\a_1=0$
(we use freedom of renumbering the spectral values).

The conjecture from the previous section is equivalent to the claim
that other $\a_i$ equal either $0$ or $\infty$. But this will be
handled in a separate paper.

\section{Integrating the system along characteristics}

 \abz
Let us now consider the quotient of the submanifold $\{H=0\}\subset
J^1(\R^n)$ by the Cauchy characteristics. We can identify it with
the transversal section $\Sigma^{2n-1}=\{H=0,u=\op{const}\}$. The
solutions will be $(n-1)$-dimensional manifolds of the induced
exterior differential system.

Note that we should augment the system with the contact form
$\oo=du-\sum p_i\,dx^i$ and its differential $\Omega_0=\sum dx^i\we
dp_i$. Note that if we choose $f(u)$ to be the solution of the ODE
(\ref{ODE}), then $\frac1{n!}\Omega_n=dp_1\we\dots\we dp_n=0$ on
solutions.

Let us start investigation from the case $\mathbf{n=2}$. In this
case the induced differential system is given by two 1-forms:
 $$
\theta=i_{X_H}\Omega_1|_\Sigma=p_1\,dp_2-p_2\,dp_1-\frac1u(p_1\,dx_2-p_2\,dx_1)
 $$
and $\theta_0=i_{X_H}\Omega_0|_\Sigma=p_1\,dp_1+p_2\,dp_2$, but it
vanishes on $\Sigma$. The form $\theta$ is contact: $\theta\we
d\theta\ne0$, so solutions of $\E$ are represented by all Legendrian
curves on $(\Sigma^3,\theta)$.

Consider now $\mathbf{n=3}$. In this case we know that
$\op{Sp}(A)=\{0,\frac1{u-\a},\frac1{u+\a}\}$ (in fact, $\a=0$, but
let us pretend we do not know it yet).

We have: $f=\frac{2u}{u^2-\a^2}$, $f'+f^2=\frac{2}{u^2-\a^2}$.

Again $\theta_0=i_{X_H}\Omega_0$ vanishes on $\Sigma^5$, so the
exteriour differential system is generated by two 2-forms:
 \begin{gather*}
\theta_1=i_{X_H}\Omega_1=(p_1\,dp_2-p_2\,dp_1)\we dx_3+ (p_2\,dp_3-p_3\,dp_2)\we dx_1\\
+(p_3\,dp_1-p_1\,dp_3)\we dx_2-\tfrac{2u}{u^2-\a^2}(p_1\,dx_2\we dx_3 +p_2\,dx_3\we dx_1 +p_3\,dx_1\we dx_2);\\
\theta_2=\tfrac12i_{X_H}\Omega_2=p_1\,dp_2\we dp_3 +p_2\,dp_3\we
dp_1 +p_3\,dp_1\we dp_2
\qquad\\
\qquad -\tfrac1{u^2-\a^2}(p_1\,dx_2\we dx_3 +p_2\,dx_3\we dx_1
+p_3\,dx_1\we dx_2).
 \end{gather*}

The integral surfaces of this system integrate to solutions of $\E$.

 {\bf Digression.}
Let us choose another section for $\Sigma'\subset J^1(\R^3)$: since
the Cauchy characteristics are given by the system $\{\dot
x_i=p_i,\dot u=1\}$, we can take in the domain $p_3>0$:
$x_3=\op{const}$, $p_3=\sqrt{1-p_1^2-p_2^2}$. Then the forms giving
the differential system are given by (being multiplied by $p_3$):
 \begin{gather*}
\theta'_1= \bigl((1-p_2^2)\,dp_1+p_1p_2\,dp_2\bigr)\we dx_2+ dx_1\we\bigl(p_1p_2\,dp_1+(1-p_1^2)\,dp_2\bigr)\\
\qquad\qquad -\tfrac{2u}{u^2-\a^2}(1-p_1^2-p_2^2)\,dx_1\we dx_2;\\
\theta'_2= dp_1\we dp_2-\frac{1-p_1^2-p_2^2}{u^2-\a^2}\,dx_1\we
dx_2.\qquad\qquad
 \end{gather*}
If we identify $\Sigma'\simeq J^1(\R^2)$ with the contact form
$\oo'=du-p_1dx^1-p_2dx^2$, the above 2-forms become represented by
the following Monge-Ampere equations:
 \begin{gather*}
(1-u_x^2)u_{xx}+2u_xu_y\cdot u_{xy}+(1-u_x^2)u_{yy}=\tfrac{2u}{u^2-\a^2}(1-u_x^2-u_y^2),\\
u_{xx}u_{yy}-u_{xy}^2=\tfrac1{u^2-\a^2}(1-u_x^2-u_y^2).
 \end{gather*}
Compatibility of this pair yields $\a=0$.

 \begin{rk}
The above system is of the kind investigated in \cite{KL$_3$}: when
the surface $\Sigma^2=\op{graph}\{u:\R^2\to\R^1\}\subset\R^3$ has
prescribed Gaussian and mean curvatures, $K$ and $H$ respectively
(this leads to a complicated overdetermined system). In fact the
PDEs of the above system can be written in the form $H=F_1(u,\nabla
u),K=F_2(u,\nabla u)$.
 \end{rk}


 \vspace{-10pt} \hspace{-20pt} {\hbox to 12cm{ \hrulefill }}
\vspace{-1pt}

{\footnotesize \hspace{-10pt} Institute of Mathematics and
Statistics, University of Troms\o, Troms\o\ 90-37, Norway.

\hspace{-10pt} E-mails: \quad kruglikov\verb"@"math.uit.no, \quad
lychagin\verb"@"math.uit.no.} \vspace{-1pt}

\end{document}